\documentclass[11pt] {amsart}
\usepackage{amsmath,amsfonts,amssymb,amsxtra,latexsym,amscd,enumerate,amsthm}

\usepackage{latexsym}
\usepackage{epsfig}
\usepackage{verbatim}

\def\r{\mathcal{R}}
\def\g{\gamma}
\def\a{\alpha}

\def\d{\delta}

\def\p{\Phi}

\def\ep{\epsilon}
\def\es{\emptyset}

\def\o{\omega}

\def\R{\mathbb{R}}
\def\C{\mathcal{C}}
\def\I{\mathcal{I}}
\def\Q{\mathcal{Q}}

\def\K{\mathcal{K}}
\def\D{\mathcal{D}}
\def\P{\mathbb{P}}
\def\p{\mathcal{P}}

\def\W{\mathcal{W}}

\def\M{\mathcal{M}}
\def\E{\mathcal{E}}
\def\S{\mathcal{S}}

\def\A{\mathcal{A}}
\def\B{\mathcal{B}}
\def\TT{\mathbb{T}}
\def\N{\mathbb{N}}

\def\beq{\begin{equation}}
\def\eeq{\end{equation}}

\def\beq{\begin{equation}}
\def\eeq{\end{equation}}

\newtheorem{d0}{Definition}

\newtheorem{l1}{Lemma}
\newtheorem{p1}{Proposition}

\begin{document}
\title[]{On Stein's Conjecture on the Polynomial Carleson Operator}

\author{Victor Lie}

\date{\today}

\address{Department of Mathematics, UCLA, Los Angeles CA 90095-1555}

\email{vilie@math.ucla.edu}
\address{Institute of Mathematics of the Romanian Academy, Bucharest, RO
70700 \newline \indent  P.O. Box 1-764}

\keywords{Time-frequency analysis, Carleson's Theorem, Stein's
conjecture, polynomial phase.}

\subjclass[2000]{42A20, 42A50.}

\maketitle
\begin{abstract}
We prove that the generalized Carleson operator $C_d$ with
polynomial phase function is of strong type $(p,r)$, $1<r<p<\infty$;
this yields a positive answer in the $1<p<2$ case to a conjecture of
Stein which asserts that for $1<p<\infty$ we have that $C_d$ is of
strong type $(p,p)$. A key ingredient in this proof is the further
extension of the {\it relational} time-frequency perspective
(introduced in \cite{q}) to the general polynomial phase.
\end{abstract}
\section{\bf Introduction}

We define the (generalized) polynomial Carleson operator as
\beq\label{v3} C_{d}f(x):=\sup_{Q\in\Q_{d}}\left|
\,p.v.\int_{\TT}\frac{1}{y}\,e^{i\,Q(y)}\,f(x-y)\,dy\,\right|\:,
\eeq where here $d\in \N$, $\Q_{d}$ is the class of all real
polynomials $Q$ with $\textrm{deg}(Q)\leq d$, and $f\in C^{1}(\TT)$
($\TT=[-\frac{1}{2},\frac{1}{2}]$). Remark that in the case $d=1$,
$C_1$ corresponds to the classical Carleson operator.

The main result of this paper is:

$\newline$
{\bf Theorem.} {\it Let $1< r< p < \infty$; then
$$\left\|C_{d}f\right\|_{L^{r}(\TT)}\lesssim_{p,r,d}\left\|f\right\|_{L^{p}(\TT)}.$$}

Furthermore, combining this result with the methods in
\cite{6} and with some general interpolation techniques, we
obtain:

$\newline$
{\bf Corollary.} {\it i) If $1<p\leq 2$ then $C_d$ is of weak type $(p,p)$.

ii) If $1<p<2$ then $C_d$ is of strong type $(p,p)$.}

$\newline$ \indent As one may observe, i) extends Carleson's theorem
on the pointwise convergence of Fourier series, which asserts that
$C_1$ is of weak type $(2,2)$ (\cite{1}, \cite{2}, \cite{4}). Also,
for $1<p<2$, ii) recovers (for $d=1$) the further extension of Hunt
(\cite{3}) and in the same range of exponents gives (for general
$d$) a positive answer to:

$\newline$
{\bf Conjecture (Stein \cite{7},\cite{8}).}
{\it If $1<p<\infty$ then $C_d$ is of strong type $(p,p)$.}

$\newline$ \indent Our results are heavily based on the intuition
and methods developed in \cite{q}, which further were significantly
influenced by the powerful geometric and combinatorial ideas
presented in \cite{2}.

We recall here that one key geometric ingredient in the proof in
\cite{q} was to regard the quadratic symmetry from a {\it
relational} perspective.  As the name suggests, this perspective
stresses the importance of \emph{interactions} between objects
rather than simply treating them independently (for details, see
\cite{q}, Section 2).  This approach had as a consequence the
splitting of the operator $C_2$ into ``small pieces" with
time-frequency portrait (morally) localized near parallelograms
(tiles) of area one. In this article, following the above-mentioned
perspective, our tiles (that will reflect the time-frequency
localization of the ``small pieces" of $C_d$) will be some ``curved
regions" representing neighborhoods\footnote{For the exact meaning
of this description, see Section 2.} of polynomials in the class
$\Q_{d-1}$.

In Section 2 we present the notations and the general procedure of
constructing our tiles, in Section 3 we elaborate on the
discretization of our operator $C_d$, while Section 4 is dedicated
to the study of the interaction between tiles. The main ingredients
that are presented in Section 5 - their proofs will be postponed to
Section 7 - will help us to prove our Theorem in Section 6. In the
Appendix we include several useful results regarding the
distribution and growth of polynomials.

Finally, we mention that this paper is closely connected to
\cite{q}. Indeed, its entire conception and realization  refer to
and require knowledge from the latter. Following this, we will
maintain the same structure of the presentation and stress only the
sensitive points that differ and/or bring some further significant
new insight into the treated topic.

{\bf Acknowledgements.} I would like to thank Terence Tao and
Christoph Thiele for reading the manuscript and giving useful
feedback.

\section{\bf Notations and construction of the tiles}

As mentioned in the introduction, we denote by $\Q_{d}$ the class of
all real polynomials of degree smaller than or equal to $d$. If not
specified, $q$ will always designate an element of $\Q_{d-1}$, while
$Q$ will refer to an element of $\Q_{d}$. When
appearing together in a proof $q$ will designate the derivative of
$Q$.

Take now the canonical dyadic grid in $[0,1]=\TT$ \footnote{The
reader should not be confused by the fact that, depending on our
convenience, the symbol $\TT$ may refer to a different unit interval
from that mentioned in the statement of our Theorem.} and in $\R$.
Throughout the paper the letters $I,\:J$ will refer to dyadic
intervals corresponding to the grid in $\TT$ while $\a^1,\ldots
\a^d$ will represent dyadic intervals associated with the grid in
$\R$. Now, if $I$ is any (dyadic) interval we denote by $c(I)$ the
center of $I$. Let $I_r$ be the ``right brother" of I, with
$c(I_r)=c(I)+|I|$ and $|I_r|=|I|$; similarly, the ``left brother" of
$I$ will be denoted $I_l$ with $c(I_l)=c(I)-|I|$ and $|I_l|=|I|$. If
$a>0$ is some real number, by $aI$ we mean the interval with the
same center $c(I)$ and with length $|aI|=a|I|$; the same conventions
apply to intervals $\{\a^k\}_k$.

A \emph{tile} $P$ is a $(d+1)$-tuple of dyadic (half open)
intervals, {\it i.e.} $P=[\a^1,\a^2,\ldots \a^d,I]$, such that
$|\a^j|=|I|^{-1}$, $j\in\{1,\ldots d\}$. The collection of all tiles
$P$ will be denoted by $\P$.

Now, for each tile $P=[\a^1,\a^2,\ldots \a^d,I]$ we will associate a
geometric (time-frequency) representation, denoted $\hat{P}$. The
exact procedure is as follows: for $I$ as before we first set
$x_{I}=(x_I^1, x_I^2,\ldots x^d_{I})\in \TT^d$ to be the $d-$tuple
defined as follows: $x_I^1, x_I^2$ are the endpoints of the interval $I$,
$x_I^3=\frac{x_I^1+ x_I^2}{2}$, $x_I^4=\frac{x_I^1+ x_I^3}{2}$, then
$x_I^5=\frac{x_I^3+ x_I^2}{2}$ and inductively (in the obvious
manner) we continue this procedure until we reach the $d$-th
coordinate. Then, define
$$\Q_{d-1}(P)=\{q\in \Q_{d-1}\:|\:q(x^j_{I})\in \a^j\:\:\forall\:j\in\{1,\ldots
d\} \}\:.$$

We will say that $q\in P$ iff $q\in \Q_{d-1}(P)$.

Finally, we set \beq\label{gtile} \hat{P}=\{(x,q(x))\:|\:x\in I\:
\&\: q\in P \}\:. \eeq

\noindent The collection of all geometric tiles $\hat{P}$ will be
denoted with $\hat{\P}$.

In the following we will also work with dilates of our tiles: for
$a>0$ and $P=[\a^1,\a^2,\ldots \a^d,I]$ we have
$aP:=[a\a^1,a\a^2,\ldots a\a^d,I]$. Similarly, we write
$$a\widehat{P}:=\widehat{aP}=\{(x,q(x))\:|\:x\in I\: \&\: q\in\Q_{d-1}(aP) \}\:.$$
Also, if $\p\subseteq \P$ then by convention
$a\p:=\left\{aP\:|\:P\in\p\right\}$; similarly, if
$\hat{\p}\subseteq \hat{\P}$ then
$\widehat{a\p}:=\left\{a\hat{P}\:|\:\hat{P}\in\hat{\p}\right\}$.

For each tile $P=[\a^1,\a^2,\ldots \a^d,I]\in\P$ we associate the
``central polynomial" $q_P\in\Q_{d-1}$ given by the Lagrange
interpolation polynomial: \beq\label{Lagrange}
 q_P(y):=\sum_{j=1}^{d}\frac{\prod_{k=1\atop{k\not=j}}^{d}(y-x_I^k)}{\prod_{k=1\atop{k\not=j}}^{d}(x_I^j-x_I^k)}\:c(\a^j)\:.
\eeq

For $P=[\a^1,\a^2,\ldots \a^d,I]$ we denote the collection of its
neighbors by $N(P)=\{P'=[{\a^1}',{\a^2}',\ldots
{\a^d}',I]\:|\:{\a^k}'\in\{\a^k,\:\a^k_r,\:\a^k_l\}\:\:\:\forall\:k\in\{1,\ldots
d\}\}$.

For any dyadic interval $I\subseteq [0,1]$, define the (non-dyadic)
intervals
$$I^{*}_{r}=[c(I)+\frac{7}{2}|I|,\,c(I)+\frac{11}{2}|I|)\:\:\:\&\:\:\:  I^{*}_l=[c(I)-\frac{11}{2}|I|,\,c(I)-\frac{7}{2}|I|)$$
and set $I^{*}=I^{*}_{r}\cup I^{*}_{l}$  and  $\tilde{I}=13I$.

Similarly, for $P=[\a^1,\a^2,\ldots \a^d,I]$ we define
$$\widehat{P^{*}_{r}}=\{(x,q(x))\:|\:x\in I^{*}_{r}\: \&\: q\in P\}$$
and further repeat the same procedure for $\widehat{P^{*}_{l}}$,
$\widehat{P^{*}}$ and $\widehat{\tilde{P}}$.

Throughout the paper $p$ will be the index of the Lebesgue
space $L^p$ and, unless otherwise mentioned, will obey $1<p<\infty$. Also,
$p'$ will be its H\"older conjugate ({\it i.e.}
$\frac{1}{p}+\frac{1}{p'}=1$), while $p^*=^{def}\min(p,p')$.

For $f\in L^p(\TT)$, we denote by
$$Mf(x)=\sup_{x\in I}\frac{1}{|I|}\int_{I}|f|$$ the Hardy-Littlewood
maximal function associated to $f$.

If $\left\{I_j\right\}$ is a collection of pairwise disjoint
intervals in $[0,1]$ and $\left\{E_j\right\}$ a collection of sets
such that for a fixed $\d\in (0,1)$ \beq \label{sets}
 E_j\subset I_j\:\:\:\:\:\&\:\:\:\:\frac{|E_j|}{|I_j|}\leq\d\:\:\:\:\:\forall\:\:j\in
 \N\:,
\eeq then we denote

\beq \label{fmax} M_{\d}f(x):=\left\{
                        \begin{array}{rl}
                        \sup_{I\supset I_j}\frac{1}{|I|}\int_{I}|f|, \  \mbox{if} \  x\in E_j \\
                        0 \qquad, \  \mbox{if} \ x \notin E_j
                        \end{array} \right.
\:.\eeq

The symbol $c(d)$ will designate a positive constant depending only on $d$; furthermore this constant is allowed to change from line to line.

For $A,\:B>0 $ we say $A\lesssim B$ (resp. $A\gtrsim B$) if there exist an
absolute constant $C>0$ such that $A<CB$ (resp. $A>CB$); if the constant $C$
depends on some quantity $\d>0$ then we may write $A\lesssim_{\d}B$.
If $C^{-1}A<B<CA$ for some (positive) absolute constant $C$
then we write $A\approx B$.

As in \cite{q}, for $x\in \R$ we set $\left\lceil
x\right\rceil:=\frac{1}{1+|x|}$.

The exponents $\eta$ and $\ep $  may change throughout the paper.

In what follows, for notational simplicity we will refer to the operator $C_d$ as $T$.

\section{\bf Discretization}

We first express $T$ in terms of its elementary building blocks:
 $$Tf(x)=\sup_{a_1,\ldots a_d \in \R}|M_{1,a_1}\ldots M_{d,a_d} H{M^{*}_{1,a_1}}\ldots {M^{*}_{d,a_d}}f(x)|=\sup_{Q\in\Q_d}|T_{Q}f(x)|\:,$$
where $\{M_{j,a_j}\}_{j\in \{1,\ldots d\}}$ is the family of
(generalized) modulations given by
$$M_{j,a_j}f(x):=e^{i a_j x^j}\:f(x)\:\:\:\:\:\:\:\:\:\:\:j\in \{1\ldots d\}$$
(here $f\in L^p,\:a_j\in\R\:\:\&\:\:x\in \TT$) and
 $$T_{Q}f(x)=\int_{\TT}{\frac{1}{y}\,e^{i\,(Q(x)-Q(x-y))}\,f(x-y)\,dy}$$ with
 $Q\in\Q_d$ given by $Q(y)=\sum_{j=1}^{d}a_j\: y^j$.

 Equivalently
$$T_{Q}f(x)=\int_{\TT}{\frac{1}{x-y}\,e^{i\,(\int_{y}^{x}q)}\:f(y)\,dy}\:,$$ where
here, as mentioned in the previous section, $q$ stands for the derivative of $Q$.

Now linearizing $T$ we write
  $$Tf(x)=T_{Q_x}f(x)=\int_{\TT}{\frac{1}{x-y}\,e^{i\,(\int_{y}^{x}q_x)}\:f(y)\,dy}\:,$$
where $Q_x(y):=\sum_{j=1}^{d}a_j(x)\: y^j$ with
$\{a_j(\cdot)\}_{j\in \{1,\dots d\}}$ measurable functions and $q_x$ refers to the derivative
of $Q_x$ ({\it i.e.} $q_x(t)=\frac{d}{dt}Q_x(t)$).

Further, proceeding as in \cite{2} and \cite{q}, we define $\psi$ to
be an odd $C^{\infty}$ function such that
$\operatorname{supp}\:\psi\subseteq \left\{y\in
\R\:|\:2<|y|<8\right\}$ and
$$\frac{1}{y}=\sum_{k\geq 0} \psi_k(y)\:\:\:\:\:\:\:\:\:\forall\:\:0<|y|<1\:,$$
where by definition $\psi_k(y):=2^{k}\psi(2^{k}y)$ (with $k\in \N$).
Using this, we deduce that
$$Tf(x)=\sum_{k\geq 0}T_{k}f(x):=\sum_{k\geq 0}\int_{\TT}e^{i\,(\int_{y}^{x}q_x)}\,\psi_{k}(x-y)\,f(y)\,dy\:.$$

Now for each $P=[\a^1,\a^2,\ldots \a^d,I]\in\P$ let
$E(P)=\left\{x\in I\:|\:q_x\in P\right\}$. Also, if $|I|=2^{-k}$ ($k\geq0$), we
define the operators $ T_P$ on $L^2(\TT)$ by
$$T_{P}f(x)=\left\{\int_{\TT}e^{i\,(\int_{y}^{x}q_x)}\,\psi_{k}(x-y)\,f(y)\,dy\right\}\chi_{E(P)}(x)\:.$$

As expected, if $\P_k:=\left\{P=[\a^1,\a^2,\ldots
\a^d,I]\in\P\:|\:|I|=2^{-k}\right\}$, for fixed $k$ the
$\left\{E(P)\right\}$ form a partition of $[0,1]$, and so
$$T_{k}f(x)=\sum_{P\in\P_{k}}T_{P}f(x)\:.$$

Consequently, we have
$$Tf(x)=\sum_{k\geq 0}T_{k}f(x)=\sum_{P\in\P}T_{P}f(x)\:.$$

This ends our decomposition.

Finally, note that (as in \cite{q}) we may assume that
$$\operatorname{supp}\:\psi\subseteq \left\{y\in\R\:|\:4<|y|<5\right\}\:.$$
Consequently, for a tile $P=[\a^1,\a^2,\ldots \a^d,I]$, the
associated operator has the properties
$$\operatorname{supp}\:T_P\subseteq I\:\:\:\:\:\:\:\:\:\:\:\&\:\:\:\:\:\:\operatorname{supp}\:T_P^{*}\subseteq \left\{x\:|\:3|I|\leq dist(x,I)\leq 5|I|\right\}=I^{*}\:,$$
where here $T_P^{*}$ denotes the adjoint of $T_P$.

Also, in what follows, (splitting
$\P=\bigcup_{j=0}^{D-1}\bigcup_{k\geq0}\P_{kD+j}$ where $D$ is the
smallest integer larger than $2d\log_{2}(2d)$) we can suppose that
if $\newline P_j=[\a^1_j,\a^2_j,\ldots \a^d_j,I_j]\in\P\:$ with
$\:j\in\left\{1,2\right\}$ such that
$|I_1|\not=|I_2|,\:\:\operatorname{then}\:\:|I_1|\leq 2^{-D}\:|I_2|$
or $|I_2|\leq 2^{-D}\:|I_1|$.

\section{\bf Quantifying the interactions between tiles}

In this section we will focus on the behavior of the expression
 \beq\label{IBT}
\left|\left\langle
T^{*}_{P_1}\:f,T^{*}_{P_2}\:g\right\rangle\right|\:. \eeq Before
this, we will need to introduce some quantitative concepts that are
adapted to the information offered by the localization of
$\{T_{P_j}\}$.

$\\\newline$ {\bf 4.1. Properties of $T_{P}$ and $T_{P}^{*}$}
$\newline$

 For $P=[\a^1,\a^2,\ldots \a^d,I]\in\P$ with $|I|=2^{-k},\:k\in
\N$, we have \beq \label{v9}
\begin{array}{rl}
        &T_{P}f(x)=\left\{\int_{\TT}\:e^{i\,(\int_{y}^{x}q_x)}\,\psi_{k}(x-y)\,f(y)\,dy\right\}\chi_{E(P)}(x)\:,  \\
    &T_{P}^{*}f(x)=\int_{\TT}\:e^{-i\,(\int_{x}^{y}q_y)}\,\psi_{k}(y-x)\,\left(\chi_{E(P)}f\right)(y)\,dy\:.
\end{array}
\eeq

Based on the {\it relational} approach developed in \cite{q} we
have: \beq\label{loc}\begin{array}{rl} &\textrm{- the time-frequency
localization of $T_{P}$ is ``morally" given by the tile $\hat{P}$;}
\\ &\textrm{- the time-frequency localization of $T_{P}^{*}$ is
``morally" given by the (bi)tile $\widehat{P^{*}}$.}
\end{array}\eeq
(Remark that, due to Lemma C (see the Appendix), one may think of
$\hat{P}$ as the $|I|^{-1}/2$ neighborhood of the graph of the
``central polynomial" $q_P$ restricted to the interval $I$.)

 $\\\newline$ {\bf 4.2. Factors of a tile} $\newline$

 For a tile $P=[\a^1,\a^2,\ldots \a^d,I]$ we define two quantities:

$\newline$a)$\:\:\:\:\:$     an {\it absolute} one (which may be
regarded as a self-interaction); we define the {\bf density
(analytic) factor of $P$} to be the expression \beq\label{m}
 A_{0}(P):=\frac{|E(P)|}{|I|}\:.
 \eeq

Notice that $A_{0}(P)$ determines the $L^{2}$ norm of $T_{P}$.

$\newline$ b)$\:\:\:\:\:\:$ a {\it relative} one (interaction of $P$
$(\hat{P})$ with something exterior to it) which is of geometric
type.

Suppose first that we are given $q\in \Q_{d-1}$ and $J$ a dyadic
interval; we introduce the quantity

$$\Delta_q(J):=\frac{\operatorname{dist}^{J}(q,0)}{|J|^{-1}}\:,$$
where we used the notations (for $q_{1},q_{2}\in\Q_{d-1}$):
$$\operatorname{dist}^{A}(q_{1},q_{2})=\sup_{y\in A}\left\{\operatorname{dist}_{y}(q_{1},q_{2})\right\}\:\:\:\:\&\:\:\:\:\operatorname{dist}_{y}(q_{1},q_{2})=\left|q_{1}(y)-q_{2}(y)\right|\:.$$

Now we define the {\bf geometric factor of $P$ $(\hat{P})$ with
respect to $q$} to be the term

$$\left\lceil \Delta_q(P)\right\rceil \big(=\frac{1}{1+ | \Delta_q(P)|}\big)\:,$$
where \beq\label{v1} \Delta_q(P):=\inf_{q_1 \in
P}\Delta_{q-q_1}(I_P)\:. \eeq

$\\\newline${\bf 4.3. The resulting estimates} $\newline$

We conclude this section by observing how the above quantities
relate in controlling the interaction in \eqref{IBT}.

As expected, we need to quantify the relative position of
$\widehat{P_1^{*}}$ with respect to $\widehat{P_2^{*}}$. (We
consider only the nontrivial case $I_{P_1}^{*}\cap I_{P_2}^*
\not=\es $; also, throughout this section we suppose that $|I_1|\geq
|I_2|$.)

\begin{d0}\label{fact} Given two tiles $P_1$ and $P_2$, we define the {\bf geometric factor of the pair ($P_1,P_2$)} by $$\left\lceil\Delta(P_1,P_2) \right\rceil\:,$$ where
$$\Delta(P_1,P_2)=\Delta_{1,2}:=\frac{\sup_{y\in I_2}\{\inf_{{q_1\in P_1}\atop{{q_2\in P_2}}}\operatorname{dist}_y(q_{1},q_{2})\}}{|I_2|^{-1}}\:.$$
 \end{d0}

With these notations, remark (using the results in the Appendix)
that we have
$${\left\lceil\Delta_{1,2}\right\rceil}\approx_{d}\max\left\{{\left\lceil\Delta_{q_{P_1}}(P_2)\right\rceil},\:{\left\lceil\Delta_{q_{P_2}}(P_1)\right\rceil}\right\}\:.$$

For $P_1$ and $P_2$ as above, we define the ``interaction
polynomial"
$$q_{1,2}:=q_{P_1}-q_{P_2}\:.$$

Fix now an interval (not necessarily dyadic) $\bar{J}\subseteq\TT$,
a polynomial $q\in \Q_{d-1}$ and three positive constants
$\eta,\:v,\: w$. In what follows we will present a general procedure
for constructing two types of critical sets associated with
$\bar{J}$, $q$, $\eta$, $v$ and $w$, denoted
$\I_s(\eta,v,q,\bar{J})$ and respectively $\I_c(\eta,w,q,\bar{J})$.

 Suppose for the moment that $q\not\in\Q_{0}$; let $J$ be the
 largest\footnote{If there are two such intervals just pick either of them.}
 dyadic interval contained in $\bar{J}$. We define
$$\M_q^{\eta}(\bar{J})=\{x\in \bar{J}\:|\:x \textrm{ is a local minimum for } |q|\:\:\&\:\:|q|(x)<\eta \}\:. $$
Now since $q\in \Q_{d-1}\setminus\Q_{0}$ we have that
$\M_q^{\eta}(\bar{J})$ is a finite set of the form
$\M_q^{\eta}(\bar{J})=\{x_j\}_{j\in \{1,\ldots r\}}$ with $r\leq
2d$. Without loss of generality we suppose that the $x_j$ are
arranged in increasing order.

Further, for each $j\in \{1,\ldots r\}$ define $\tilde{q}_j(x):=q(x)-q(x_j)$ and construct the dyadic intervals $I_j^1$,
$I_j^2$ and $I_j^3$ as follows: $I_j^1$ is the smallest dyadic interval $I$ for which $x_j\in I$ and $\Delta_{\tilde{q}_j}(I)>c(d)\, v$; $I_j^2$ is the smallest dyadic interval $I$ having the left end point equal with the right end point of $I_j^1$ and for which $\Delta_{\tilde{q}_j}(I)>c(d)\, v$; similarly, $I_j^3$ is the smallest dyadic interval $I$
having the right end point equal with the left end point of $I_j^1$ and for which $\Delta_{\tilde{q}_j}(I)>c(d)\, v$.

Now set $$\S_{j}(\eta,v,q,\bar{J})=\bigcup_{k=1}^3 I_j^k$$ and define
$$\S(\eta,v,q,\bar{J})=\bigcup_{j=1}^r
\S_{j}(\eta,v,q,\bar{J}).$$

 Also set $$\C_{j}(\eta,w,q,\bar{J})= [x_j-w, x_j+w]\cap\bar{J}$$ and further take
$$\C(\eta,w,q,\bar{J})=\bigcup_{j=1}^r \C_{j}(\eta,w,q,\bar{J}).$$

We now need to do one more step before ending our construction;
suppose that $A\subseteq \bar{J}$ is a finite union of (closed)
intervals: $A=\bigcup_{j=1}^l A_j$ with $l\in\N$, $A_j=[u_j,v_j]$
(pairwise disjoint) and $\{u_j\}$ monotone increasing. Then, setting
$A_0=A_{l+1}=\emptyset$ we define

$$\E(\bar{J},A)=\left( \bigcup_{j=1}^l A_{j}\right)\cup\left(\bigcup_{j\in\{1\ldots l+1\}\atop{|C_{j}|<|A_{j-1}|\textrm{ or }|C_{j}|<|A_{j}|}}C_j\right)\:, $$
where here the intervals $C_j$ obey the partition condition
$$\bar{J}=\bigcup_{j=1}^{l+1} \:\left( A_j\cup C_j\right)\:.$$

Finally if $q\not\in\Q_0$ then define

$$\I_s(\eta,v,q,\bar{J})=\E(\bar{J},\S(\eta,v,q,\bar{J}))$$
and
$$\I_c(\eta,w,q,\bar{J})=\E(\bar{J},\C(\eta,w,q,\bar{J})).$$

Otherwise, if $q\in\Q_0$, just set
$\I_s(\eta,v,q,\bar{J})=\I_c(\eta,w,q,\bar{J})=\emptyset$.

Fix $\ep_0\in (0,1)$. Set
$w(\bar{J})=c(d)\,|J|\,{\left\lceil\Delta_{q}(J)\right\rceil}^{\frac{1}{d}-\ep_0}$,
$v(\bar{J})=c(d)\,{\left\lceil\Delta_{q}(J)\right\rceil}^{-2\ep_0}$
and $\eta(\bar{J})=c(d)\,v(\bar{J})\,w(\bar{J})^{-1}$.

Then (using the results in the Appendix) we deduce
 \beq\label{rel}
\I_s(\eta(\bar{J}),v(\bar{J}),q,\bar{J})\subseteq\I_c(\eta(\bar{J}),w(\bar{J}),q,\bar{J}).
\eeq

 We now define the ($\ep_0$-){\bf critical intersection set}
$I_{1,2}$ of the pair $(P_1,P_2)$ as
 \beq\label{criticset}
I_{1,2}:=\I_c(\eta_{1,2},w_{1,2},q_{1,2},\tilde{I_1}\cap
\tilde{I_2})\:,\eeq where $\eta_{1,2}:=\eta(\tilde{I_1}\cap
\tilde{I_2})$ and $w_{1,2}:=w(\tilde{I_1}\cap \tilde{I_2})$.

Notice that, based on \eqref{rel} and Lemma C of the Appendix, we
have that \beq\label{ctbs}
 \bigcup_{{q_j\in P_j}\atop{j\in\{1,2\}}}\left\{y\in\tilde{I}_2\:\big|\:\frac{|q_1(y)-q_2(y)|}{|I_2|^{-1}}\leq
{\left\lceil\Delta_{1,2}\right\rceil}^{-\frac{1}{d}-\ep_0}\right\}\subseteq
I_{1,2}\:.\eeq

Now using \eqref{ctbs} together with the principle of
(non-)stationary phase, one deduces the following:

$\newline$ {\bf Lemma 0.} {\it Let $P_1\:,\:P_2\:\in\P$; then we
have \beq\label{v15} \left|\int
\tilde{\chi}_{I_{1,2}^c}T_{P_1}^{*}f\:\overline{T_{P_2}^{*}g}\:\right|\lesssim_{\:n,\:d,\:\ep_0}{\left\lceil
\Delta(P_1,P_2)\right\rceil}^{n}\:\frac{\int_{E(P_1)}
|f|\int_{E(P_2)}|g|}{\max\left(|I_1|,|I_2|\right)}\:\:\:\:\:\:\forall\:n\in
\N\:, \eeq
 \beq\label{v16} \left|\int_{I_{1,2}}
T_{P_1}^{*}f\:\overline{T_{P_2}^{*}g}\:\right|\lesssim_d{\left\lceil
\Delta(P_1,P_2)\right\rceil}^{\frac{1}{d}-\ep_0}\:\frac{\int_{E(P_1)}
|f|\int_{E(P_2)}|g|}{\max\left(|I_1|,|I_2|\right)}\:, \eeq where
$\tilde{\chi}_{I_{1,2}^c}$ is a smooth variant of the corresponding
cut-off.

Applying the same methods for the limiting case $\ep_0=0$, we obtain
\beq\label{v17}
\left\|T_{P_1}{T}_{P_2}^{*}\right\|_2^{2}\lesssim_{d}\min\left\{\frac{|I_2|}{|I_1|},
\frac{|I_1|}{|I_2|} \right\}{\left\lceil
\Delta(P_1,P_2)\right\rceil}^{\frac{2}{d}}\:A_0(P_1)\:A_0(P_2)\:.
\eeq}

\begin{proof}
We first notice that relation \eqref{v16} is straightforward;
indeed, to see this we just use the relation
$$|T_{P_j}^{*}f|\lesssim
\frac{\int_{E(P_j)}|f|}{|I_j|}\,\chi_{I_j^*}\:\:\:\:\:\:\:\forall\:\:j\in\{1,2\}$$
together with the definition of $I_{1,2}$.

We turn now our attention towards \eqref{v15}. First, for notational
convenience we set $\varphi=\tilde{\chi}_{I_{1,2}^c}$; with this, we
have:
\begin{align*}
\int \varphi\: T_{P_1}^{*}f\:\overline{T_{P_2}^{*}g} & =\int f\:
\overline{T_{P_1}(\varphi\, T_{P_2}^{*}g)}\\
& =\int\int (f\chi_{E(P_1)})(x)\:
(\bar{g}\chi_{E(P_2)})(s)\:\K(x,s)\,dx\, ds\:,
\end{align*}

where
$$\K(x,s)=\int e^{i\,
[\int_{y}^{s}q_s-\int_{y}^{x}q_x]}\;\psi_{k_1}(x-y)\:\varphi(y)\:\psi_{k_2}(y-s)\;dy\:.$$

(Here we use the conventions $|I_1|=2^{-k_1}$, $|I_2|=2^{-k_2}$ with
$k_2\geq k_1$ positive integers.)

 Now making the change of variables $y=|I_2|\,u$ and using the way in
 which we defined $I_{1,2}$ and $\varphi$, we deduce that

$$\left|\K(x,s)\right|\lesssim |I_1|^{-1}\left|\int_{\TT}e^{i \phi(u)}\,r(u)\,du\right|$$
 with $r\in C_{0}^{\infty}(\R)$ such that $|\partial^{l}r(u)|\lesssim_d (\left\lceil
\Delta(P_1,P_2)\right\rceil^{\ep_0-\frac{1}{d}})^l$ ($l\in\N$) and
$\left\|\,\partial\phi\,\right\|_{L^{\infty}(\textrm{supp}\:
r)}\gtrsim_d \left\lceil
\Delta(P_1,P_2)\right\rceil^{-\ep_0-\frac{1}{d}}$. Using the
non-stationary phase principle we thus obtain \eqref{v15}.

For $\eqref{v17}$, we set $\ep_0=0$ in the previous argument.

\end{proof}

\section{\bf Main ingredients}

In this section, we will present the concepts and results that we
need for proving our theorem.

 \begin{d0}\label{mass}                     For $P=[\a^1,\a^2,\ldots \a^d,I]\in\P$ we define the {\bf mass} of $P$
 to be

\beq\label{v18} A(P):=\sup_{{P'=[{\a^1}',{\a^2}',\ldots
{a^d}',I']\in\:\P}\atop{I\subseteq
I'}}\frac{|E(P')|}{|I'|}\:\left\lceil
\Delta(2P,\:2P')\right\rceil^{N}\:, \eeq where $N$ is a fixed large
natural number.
\end{d0}

Next, we introduce a qualitative concept that characterizes the
overlapping relation between tiles.
 \begin{d0}\label{ord}
   Let $P_j=[\a^1_j,\a^2_j,\ldots \a^d_j,I_j]\in\P$ with $j\in\left\{1,2\right\}$. We say that
$\newline$ - $P_1\leq P_2$       iff       $\:\:\:I_1\subseteq I_2$
and      $\exists\:\:q\in P_2\:\:s.t.\:\:q\in P_1\:,$ $\newline$ -
$P_1\trianglelefteq P_2$     iff     $\:\:\:I_1\subseteq I_2$  and
$\forall \:\:q\in P_2\:\:\Rightarrow\:\:q\in P_1\:.$
\end{d0}

 \begin{d0}\label{tree}

 We say that a set of tiles $\p\subset\P$ is a \emph{tree} (relative to $\leq$) with top $P_0$ if
the following conditions are satisfied: $\newline
1)\:\:\:\:\:\forall\:\:P\in\p\:\:\:\Rightarrow\:\:\:\:\frac{3}{2}P\leq
10 P_0 $ $\newline 2)\:\:\:\:\:$if $P\in\p$ and $P'\in N(P)$ such
that $\frac{3}{2}P'\leq P_{0}$ then $P'\in\p$ $\newline
3)\:\:\:\:\:$if $P_1,\:P_2\: \in\p$ and $P_1\leq P \leq P_2$ then
$P\in\p\:.$
\end{d0}

Now we can state the main results of this section; their proofs will
be postponed until Section 7.

\begin{p1}\label{prop1}
There exists $\eta\in(0,1/2)\:$ (depending only on the degree $d$)
s.t. if $\:\p\:$ is any given family of disjoint tiles (i.e. no two
of them can be related through $\leq$) with the property that
$$ A(P)\leq\d\:\:\:\:\:\forall\:\:\:P\in\p$$
then, for $1<p<\infty$, we have
$$\left\|T^{\p}\right\|_{p}\lesssim_{p,d}\d^{\eta(1-\frac{1}{p^*})}\:.$$
\end{p1}

\begin{p1}\label{prop2}
Let $\left\{\p_j\right\}_j$ be a family of trees with tops
$P_j=[\a^1_j,\a^2_j,\ldots \a^d_j,I_j]$. Suppose that
$\newline\:1)\:\:\:\:A(P)<\d\:\:\:\:\:\:\:\forall\:j,\:P\in\p_j\:,$
$\newline\:2)\:\:\:\forall\:\:
k\not=j\:\:\&\:\:\forall\:\:P\in\p_j\:\:\:\:\:\:\:\:2P\nleq
10P_k\:,$ $\newline\:$3)    No point of $[0,1]$ belongs to more than
$K\d^{-(1+\rho)}$ of the $I_j$ (here $\rho$ is a fixed number with
$0<\rho\leq\min\left\{1,\frac{1}{2} | \frac{p-1}{2-p}| \right\}\:$).

Then there exist a constant $\eta\in(0,\frac{1}{2})$ (depending only
on $d$) and a set $F\subset\TT$ with $|F|\lesssim\d^{50}K^{-M}$
(here $M\in\N$ is fixed) such that $\forall\:f\in L^p(\TT)$ we have
$$\left\|\sum_{j}T^{{\p}_j}f\right\|_{L^p(F^{c})}\lesssim_{p,d} \left(\d^{\eta(1-\frac{1}{p^*})}M\log{K}+K^{\frac{1}{p^*}-\frac{1}{p'}}\:\d^{\frac{1}{p}-(1+\rho)(\frac{1}{p^*}-\frac{1}{p'})}\right)\left\|f\right\|_p\:.$$

(Remark: Any collection of tiles $\p$ that can be represented as
$\cup_{j}\p_j$ with the family $\left\{\p_j\right\}$ respecting the
conditions mentioned above will be called a ``forest".)
\end{p1}

\section{\bf Proof of ``pointwise convergence"}

We now present the proof of our Theorem.

The main challenge will be to organize the collection $\P$ into a
``controlled number" of forests. For this, the first step is to
split our family of tiles into a union of subfamilies having uniform
density. More exactly, (as in \cite{2}, \cite{q}) we decompose
$$\P=\bigcup_{n=0}^{\infty}\p_{n}\:\:\:\:\:\:\:\textrm{with}\:\:\:\:\:\p_{n}=\left\{P\in\P\:|\:2^{-n-1}< A(P)\leq 2^{-n}\right\}\:.$$
In what follows, writing $T\:=\:\sum_{n=0}^{\infty}T^{\p_n}$, we
will show that
$$\exists\:\:\eta=\eta(p,d)>0\:\:\:\textrm{such that}\:\:\:\left\|T^{\p_{n}}\right\|_{L^p(E_n^c)}\lesssim_{p,d} 2^{-\eta n}\:,$$
where $E_n$ are some exceptional (``small") sets inside the torus.

 Now, for a fixed $n$, we need to search for clustered families
of tiles - the natural environments for the future trees - living
inside $\p_n$. Consequently, the next step is to identify some
preferred tiles that will help determine the tops of these trees.
For this, we select the tiles $\left\{\bar{P}_{k}\right\}$
$(\bar{P}_k=[\a^1_k,\a^2_k,\ldots \a^d_k,\bar{I}_k])$ to be those
$(d+1)$-tuples which are maximal with respect to $\leq$ and obey
$\frac{|E(P)|}{|I_P|}\geq 2^{-n-1}$. Proceeding as in \cite{q}, (and
also maintaining the same notations) we define \beq \label{triang}
 \p_n^{0}:=\left\{P\in\p_n\:|\:\:\exists\:k\in \N \:s.t.\:\:\:\:\:4P\triangleleft \bar{P}_{k}\right\}
 \eeq
and set $$\C_n:=\left\{P\in\p_n\:|\:\operatorname{there\:are\: no\:
chains}\:P\lneq P_{1}\lneq\ldots\lneq
P_{n}\:\&\:\left\{P_j\right\}_{j=1}^{n}\subseteq\p_n\:\right\}\:.$$
With this done, it is easy to see that
$$\p_n\setminus\C_n\subseteq\p_n^{0}\:.$$ Now defining the set
$\D_n\subseteq \C_n$ with the property $\p_n\setminus\D_n=\p_n^{0}$,
we remark that $\D_n$ breaks up as a disjoint union of a most $n$
sets $\D_{n1}\cup\D_{n2}\cup\ldots\cup\D_{nn}$ with no two tiles in
the same $\D_{nj}$ comparable. Applying Proposition 1, we see that
\beq \label{D}
\left\|T^{\D_{n}}\right\|_p\leq\sum_{j=1}^{n}\left\|T^{\D_{nj}}\right\|_p\lesssim_{p,d}\sum_{j=1}^{n}2^{-n\eta(1-\frac{1}{p^*})}\lesssim
2^{-n\eta(1-\frac{1}{p^*})} \:, \eeq and so, in what follows, it
will be enough to limit ourselves to the set of tiles $\p_n^{0}$.

Once we have established a certain structure of our collection of
tiles - all our elements are clustered near some maximal elements
$\left\{\bar{P}_{k}\right\}$ - to move closer towards the concept of
``forest" we need to learn how to control the number of these
maximal elements (see 3) Proposition 2). For this, we make the
following reasoning: Define the counting function
$N(x)=^{def}\sum_{k}\chi_{\bar{I}_k}(x)$ and set
$$G_n:=\left\{x\in\TT\:|\:x\:\operatorname{is\:contained\:in\:more\:than\:}2^{(1+\rho)n}K\:\operatorname{of\:the\:}|\bar{I}_{k}|
\right\}\:,$$ where here the parameter $\rho$ is as defined in the
statement of Proposition 2.

In contrast with the $L^2$ case, to obtain $L^p$ bounds for our operator $T$,
we need to be more precise in
estimating the measure of $G_n$, or equivalently to better control
the size of the level sets corresponding to $N$:
$$G_n=\{x\in\TT\:|\:N(x)\geq 2^{(1+\rho)n}K\}.$$

For this task we may proceed in two ways: the first (``modern") one
uses $BMO$- estimates and is shorter, more efficient but less
descriptive; the second (``classical") one, may be found in \cite{2}
(Section 8) and uses a vector-valued variant of the Hardy-Littlewood
maximal theorem.

{\it First approach:}

We first observe that from the definition of $\{\bar{I}_k\}$ we have
that $N\in BMO(D)$\footnote{Throughout the paper we will denote with
$BMO(D)$ the dyadic $BMO$ on the torus.} with $\|N\|_{BMO(D)}\leq
2^n$. Applying now the John-Nirenberg inequality, we have
$(\gamma>0)$
$$\left|\,\left\{x\in\TT\:\left| \:|N(x)-\int_{\TT}N|> \gamma\right\}\,\right|\right.\lesssim
e^{-\frac{c\gamma}{\|N\|_{BMO(D)}}}\:,$$ where $c>0$ is an absolute
constant.

Now taking into account the fact that $\|N\|_{1}\lesssim 2^n$ and
setting $\gamma:=2^{(1+\rho)n}K$ we deduce that
$$|G_n|\lesssim e^{-c2^{\rho n}K}\:. $$

{\it Second approach:}

Since by Chebyshev's inequality \beq\label{ceb} |\{x\:|\:N(x)\geq
\gamma\}|\leq \gamma^{-r} \int N^r\:\:\:\:\:\:\forall\:\:r>0\:, \eeq
we are led to the study of the $L^r$ bounds for $N$, where now our
interest is for $r$ to be large.

The strategy for this is to relate the behavior of
$\{\chi_{\bar{I}_k}\}$ with that of $\{\chi_{E(\bar{P}_k)}\}$, since
the latter has better disjointness properties. Indeed, set
$h_k(x):=\chi_{E(\bar{P}_k)}$ and $h^{*}_k=M(h_k)$ the
Hardy-Littlewood maximal operator applied to $h_k$; then we have
\beq\label{maxestim} N(x)\lesssim_{\rho'}
2^{n(1+\rho')}\sum_{k}(h_k^{*})^{1+\rho'}(x)\:\:\:\:\:\forall\:x\in\TT\:,\eeq
where $\rho'<\rho$ is a fixed positive constant.

Now applying the $l^p$ vector-valued variant of the Hardy-Littlewood
maximal theorem $(r>1)$, we deduce

\beq\label{VHL}
\int\left[\sum_{k}(h_k^{*})^{1+\rho'}\right]^r\lesssim_{r,\rho'}\int\left[\sum_{k}(h_k)^{1+\rho'}\right]^r=
\sum_{k}|E(\bar{P}_k)|\leq 1\:,\eeq and so, by \eqref{ceb},
\eqref{maxestim} and \eqref{VHL} we have

$$|\{x\:|\:N(x)\geq \gamma\}|\leq \gamma^{-r}2^{nr(1+\rho')}\:.$$

Now choosing $\gamma=2^{(1+\rho)n}K$ and $r=M$ we conclude that
there exists $\eta=\eta(p,d)>0$ such that
$$|G_n|\lesssim_{p} \frac{2^{-n\eta}}{K^M}\:.$$

$\newline$ Now, to control the number of tops that lie one upon the
other, we will use the estimates (of the set $G_n$) just obtained to
erase a few more tiles from $\p_n^0$. More exactly if we set
$$\p_n^{G}=\left\{P=[\a^1,\a^2,\ldots \a^d,I]\in\p_n^{0}\:|\:I\nsubseteq G_n\right\}\:,$$ we have that
\beq \label{exc}
T^{\p_n^{G}}f(x)\:=\:T^{\p_n^{0}}f(x)\:\:\:\:\:\:\:\:\:\:\forall\:f\in
L^{p}(\TT)\:\:\&\:\:x\in G_n^{c}\:. \eeq Now, since we have good
control on the measure of the set $G_n$, we will erase from
$\left\{\bar{P}_{k}\right\}$ all $\bar{P}_k$ with
$\bar{I}_k\subseteq G_n$ and focus further on estimating
$T^{\p_n^{0}}$ only on $G_n^c$.

The point is that now the set $\p_n^{G}$ has the following properties: $\newline 1)\:\:A(P)\leq
2^{-n}\:\:\:\:\:\:\:\:\forall\:\:P\in\p_n^{G}\:,$ $\newline
2)\:\:\forall\: P\in\p_n^{G}\:\:\Rightarrow\:\exists\:k\in N
\:st\:\:\:\:\:4P\trianglelefteq \bar{P}_{k}\:,$ $\newline 3)\:\:$No
$x\in\TT$ belongs to more than $2^{(1+\rho)n}K$ of the $\bar{I}_k$'s.

We remark that we are coming closer to the definition of a forest -
see Proposition 2. One ingredient that is still missing refers to
the separation condition 2). This will force us to rearrange the
tiles inside $\p_n^{G}$ into collections that are clustered near
some dilated maximal elements. The exact procedure follows the same
lines as those described in the section 7.2. of \cite{q}; for
convenience, we will sketch the major steps in the next section,
leaving the details for the reader.

$\newline${\bf 6.2. Decomposing into forests}$\newline$

As announced, in this section we intend to reorganize the set
$\p_n^{G}$; more precisely, we will show that up to a
negligible\footnote{The set can be written as a union of at most
$c(d)$ families of disjoint tiles.} set of tiles
$$\p_n^{G}\cong\bigcup_{j=1}^{M}\B_{nj}\:,$$
with each $\B_{nj}$ a forest and $M$ some constant less than $2n\log K$.
As in \cite{q} we set
$$B(P):=\#\left\{j\:|\:4P\trianglelefteq \bar{P}_j\right\}\:\:\:\:\:\:\:\forall\:\:P\in\p_n^{G}\:,$$
and
$$\p_{nj}:=\left\{P\in\p_n^{G}\:|\:2^j\leq B(P)< 2^{j+1}\right\}\:\:\:\:\:\:\:\:\:\:\forall\:j\in \left\{0,..M\right\}\:.$$
Fixing now a family $\p_{nj}$ we look for candidates for the tops of
the future trees. More exactly, take
$\left\{P^{r}\right\}_{r\in\left\{1,\ldots
s\right\}}\subseteq\p_{nj}$ to be those tiles with the property that
$\:4P^{r}$ are maximal\footnote{Here we use the following
convention: let be $\D$ a collection of tiles; $P$ is maximal
(relative to $\leq$) in $\D$ iff $\forall\:\:P'\in \D$ such that
$P\leq P'$ we also have $P'\leq P$.} elements with respect to the
relation $\leq$ inside the set $4\p_{nj}$.

Now, in all the reasonings that we will make further, we will use
the following four essential properties:
\begin{enumerate}[(A)]
    \item $\:\:\:\: 4P^{l}\leq 4P^{k}\:\:\Rightarrow\:\:I_l=I_k\:; \\$
    \item $\:\:\:\: \forall\:P\in\p_{nj}\:\:\:\exists\:\:\:P^{l}\:\:\:\textrm{s.t.}\:\:\:\:4P\leq
4P^{l}\:;\\$
  \item  $\:\:\:\:\operatorname{If}\:P\in\p_{nj}\:\:\operatorname{s.t.}\:\:\exists\:\:k\not=l\:\:\:\operatorname{with}\:\:\:\left\{{4P\trianglelefteq\:4P^{l}}\atop{4P\trianglelefteq\: 4P^{k}} \right.\:,\:\operatorname{then}\:\:\:\left\{{4P^{k}\leq4P^{l}}\:\atop{4P^{l}\leq 4P^{k}} \right.
\:;\\$
  \item  $\:\:\:\:\textrm{If}\: P_j=[\a^1_j,\a^2_j,\ldots \a^d_j,I_j]\in\P\:\:\textrm{with}
\:\:j\in\left\{1,2\right\}\: \textrm{s.t.}\:
|I_1|\not=|I_2|,\\ \operatorname{then}\:\:|I_1|\leq 2^{-D}\:|I_2|
\:\textrm{or}\: |I_2|\leq 2^{-D}\:|I_1|\:.\\$
\end{enumerate}

(While $(A),\:(B)$ and $(C)$ are derived from the definition of
$\p_{nj}$ and the way in which we have chosen
$\left\{P^{r}\right\}_{r\in\left\{1,\ldots s\right\}}$, property
$(D)$ follows from the assumption made at the end of Section 3.)

The next step, is to discard the tiles that are not ``close" to our
new maximal elements (remember that our final goal is to construct
``well separated" trees); for technical reasons, we also get rid of
the maximal elements together with their neighbors and respectively
of the minimal elements; more exactly, we define:
$$\A_{nj}:=\left\{P\in\p_{nj}\:|\:
\forall\:\:P^{l}\:\:\Rightarrow\:\:\frac{3}{2}P\nleqslant
P^{l}\right\}\cup$$
$$\left\{P\in\p_{nj}\:|\:\exists\:l\:st\:|I_P|=|I_{P^l}|\:,\:\frac{3}{2}P\leq P^{l}\right\}
\cup\left\{P\in\p_{nj}\:|\:P\:\textrm{minimal}\atop{\textrm{with respect of} ``\leq"}\right\}\:.$$

Also, we set $$\B_{nj}:=\p_{nj}\setminus\A_{nj}\:.$$
Now, using the properties (B) and (D), it is easy to see that $\A_{nj}$ forms a negligible set of tiles.

For the remaining set $\B_{nj}$, one should follow the steps bellow
(here we use the properties (A)-(D)):

 \begin{itemize}
\item  Set $S_k=\left\{P\in\B_{nj}\:|\:\frac{3}{2}P\leq P^{k}\right\};$ without loss of generality
we may suppose $\B_{nj}=\bigcup_{k=1}^{s}S_k$;
\item  Introduce the following relation among the sets $\left\{S_k\right\}_k$:
$$S_k\propto S_l$$
if $\exists\: P_1\in S_k$ and
$\exists\: P_2\in S_l$ such that $2P_1\leq 10P^l$ or $2P_2\leq 10P^k$;
\item  Deduce that $``\propto"$ becomes an order relation and that
$$S_k\propto
S_l\:\:\:\Rightarrow\:\:\:4P^k\leq4P^l\:\:\Rightarrow\:\:I^k=I^l\:;$$
\item  Let $\hat{k}:=\left\{m\:|\:S_m\propto S_k \right\}$; then
the cardinality of $\hat{k}$ is at most $c(d)$, and for
$$\hat{S}_k:=\bigcup_{m\in\hat{k}}S_m\:,$$
one has that $\hat{S}_k$ is a tree having as a top any $P^l$ with $l\in\hat{k}$.
\item  Conclude that the relation $``\propto"$, can be meaningfully extended\footnote{In this case,
the role played by the maximal element $P^k$ in the initial
definition is now taken by the top of the corresponding tree.} among
the sets $\left\{\hat{S}_k\right\}_k$ and that
$$\hat{S}_k\propto\hat{S}_l\:\:\Rightarrow\:\:k=l\:.$$
\end{itemize}
$\newline$
Consequently, from the algorithm just described, we deduce that the set $$\B_{nj}=\bigcup_{k}\hat{S}_k$$
is a forest as defined in Proposition 2.

\newpage

$\newline${\bf 6.3. Ending the proof} $\newline$

To complete the proof we proceed as follows:

We first apply Proposition 2 for each family $\B_{nj}$ and obtain that
$$\left\|T^{{\p}_{nj}}f\right\|_{L^{p}(F^{c}_{nj})}\lesssim_{p,d} 2^{-n\:\eta(p,d)}\left(M\log{K}+K^{\frac{1}{p^*}-\frac{1}{p'}}\right)\left\|f\right\|_p\:,$$
where $\eta(p,d)$ is a constant depending only on $p,\:d$, and
$F_{nj}$ is a small set with measure $|F_{nj}|\lesssim
\frac{2^{-n}}{K^M}$. As a result, denoting $F_{n}=\bigcup_{j}
F_{nj}$, we have that

\begin{align*}
\left\|T^{{\p}_n^{G}}f\right\|_{L^{p}(F^{c}_{n})} & \leq
\sum_{j=1}^{2n\log
K}\left\|T^{{\p}_{nj}}f\right\|_{L^{p}(F^{c}_{nj})} \\
 & \lesssim_{p,d} n
2^{-n\:\eta(p,d)}\left(M(\log{K})^2+K^{\frac{1}{p^*}-\frac{1}{p'}}
\log{K}\right)\left\|f\right\|_p\:,
\end{align*}
 with $|F|\lesssim \frac{n\log
K}{2^n K^M}$.

Therefore, we deduce
$$\left\|T^{{\p}_n}f\right\|_{L^{p}(E^{c}_{n})}\lesssim_{p,d} n
2^{-n\:\eta(p,d)}\left(M(\log{K})^2+K^{\frac{1}{p^*}-\frac{1}{p'}}
\log{K}\right)\left\|f\right\|_p\:,$$ where $E_{n}=F_{n}\cup G_{n}$
still has measure $\lesssim \frac{n\log K}{2^n K^M}$.

 Summing now over $n$, we obtain
$$\left\|Tf\right\|_{L^{p}(E^{c})}\lesssim_{p,d} \left(M(\log{K})^2+K^{\frac{1}{p^*}-\frac{1}{p'}}
\log{K}\right) \left\|f\right\|_p$$ with $E=\bigcup_{n}E_n$ and
$|E|\lesssim\frac{\log K}{K^M}$.

In conclusion, given $\g\:>\:0$, we have that for all $K,\:M>100$

\begin{align*}
\left|\left\{|Tf(x)|\:>\:\g\right\}\right| & \leq  \frac{\left\|Tf\right\|^{p}_{L^{p}(E^{c})}}{{\g}^p}\:+\:|E|\\
 & \lesssim_{p,d}  \left(M^p(\log{K})^{2p}+K^{\frac{p}{p^*}-\frac{p}{p'}}
(\log{K})^p\right)\frac{\left\|f\right\|_p^{p}}{{\g}^p}\:+\:\frac{\log
K}{K^M}\:.
\end{align*}

Now, if we pick $K$ to minimize the right-hand side and make $M$
large enough, we arrive at the relation
$$\left|\left\{|Tf(x)|\:>\:\g\right\}\right|\lesssim_{p,d,\ep}\left(\frac{\left\|f\right\|_p}{{\g}}\right)^{p-\ep}\:\:\:\:\:\:\:\:\:\:\:\forall\:\:\:1<p<\infty\:\:\&\:\:\ep\in\:(0,1)\:,$$
which further implies (using interpolation)
$$\left\|Tf\right\|_r\lesssim_{r,p,d} \left\|f\right\|_p\:\:\:\:\:\:\:\:\:\:\:\:\:\forall\:\:\:\:1<r<p<\infty\:,$$
ending the proof of our theorem.

\section{\bf Some technicalities - the proofs of Propositions 1 and
2}

 $\newline$ {\bf Proof of Proposition 1}
$\newline$

Following exactly the same steps\footnote{The proof is based on two
cases: for the diagonal term, since we don't have oscillatory
behavior, one uses maximal methods, while for the off-diagonal term
we use the $TT^{*}$-method together with Lemma 0.} from the
corresponding proof of Proposition 1 in \cite{q}, one can show that
there exists $\eta$ (depending only on the degree $d$) such that
$$\left\|T^{\p}\right\|_2^2\lesssim_{d} \d^{\eta}\:.$$
Now, for the case $2<p<\infty$, we observe that:
$$\left\|T^{\p}\right\|_{\infty\rightarrow\infty}=\sup_{f\in
L^{\infty}\atop{\left\|f\right\|_{\infty}\leq
1}}\left\|T^{\p}f\right\|_{\infty}\lesssim\sup_{f\in
L^{\infty}\atop{\left\|f\right\|_{\infty}\leq
1}}\sup_{P\in\p}\frac{\int_{I_{P^{*}}}|f|}{|I_P|}\lesssim 1\:,$$ so
interpolating between $L^2$ and $L^{\infty}$ we obtain the desired
conclusion.

For the case $1<p<2$ we need to focus on the behavior of
${T^{\p}}^{*}$.

\noindent Indeed, on the one hand we know that
$$\left\|{T^{\p}}^{*}\right\|_{2\rightarrow 2} = \left\|T^{\p}\right\|_{2\rightarrow
2}\lesssim_{d}\d^{\frac{\eta}{2}}\:.$$

\noindent On the other hand, for $f\in L^{\infty}$ we have
$$|{T^{\p}}^{*}f|\leq
\sum_{P\in\p}|{T^{P}}^{*}f|\lesssim\sum_{P\in\p}\frac{\int_{E(P)}|f|}{|I_P|}\chi_{I_{P^*}}\lesssim
\left\|f\right\|_{\infty}\sum_{P\in\p}\frac{|E(P)|}{|I_P|}\chi_{I_{P^*}}\:.$$

If we now add the fact
$$\left\|\sum_{P\in\p}\frac{|E(P)|}{|I_P|}\chi_{I_{P^*}}\right\|_{BMO(D)}\lesssim 1$$
we conclude that
$$\left\|{T^{\p}}^{*}\right\|_{\infty\rightarrow BMO(D)}\lesssim 1\:.$$
The claim now follows by interpolation.\footnote{We use here the
fact that $\left\|{T^{\p}}^{*}\right\|_{p'\rightarrow
p'}=\left\|{T^{\p}}\right\|_{p\rightarrow p}$.}

\begin{flushright}
$\Box$
\end{flushright}

The remainder of the section will be dedicated to proving
Proposition 2. The strategy here is as follows: for the single tree
estimate (Lemma 1) we prove the $L^p$ bounds directly; after that,
we will show only $L^2$ statements (Lemma 2, 3 and 4) that will be
just enough to prove Proposition 2 in the case $p=2$. Once we get
this, by applying interpolation techniques we will recover the
entire range of $p$.

\begin{l1}\label{tr}
Let $\delta>0$ be fixed and let $\p\subseteq\P$ be a tree such that
$$A(P)<\delta\:\:\:\:\forall\:\:\:P\in\p\:.$$ Then
\beq\label{tree}
\left\|T^\p\right\|_p\lesssim_{p,d}\delta^{\frac{1}{p}}\:. \eeq
\begin{proof}
We start by setting the parameters of our tree; more precisely, we
fix the top $P_0=[\a^1_0,\a^2_0,\ldots
\a^d_0,I_0]$, and frequency polynomial $q_0$. Since we have
specified the polynomial $q_0$ we also know the form of $Q_0$;
suppose now that
$$Q_0(y)=\sum_{j=1}^{d}a_j^0\:y^j\:.$$
Then, denoting with $g(x)={M^{*}_{1,a_1^0}}\ldots
{M^{*}_{d,a_d^0}}f(x)$ and reasoning\footnote{Here the key element
is the following perspective: ``A tree behaves like a maximal
truncated Hilbert transform".} as in the proof of Lemma 1 in
\cite{q}, one can show that:
$$|T^\p f(x)|\lesssim_d M_{\d}(R*g)(x)+M_{\d}g(x)\:,$$
where we set $R(y)=\sum_{k\in \N D}\psi_{k}(y)$ (here without loss
of generality we suppose that $\p\subset\bigcup_{k\in\N}\P_{k D}$).

Now, taking into account the fact that
$\left\|M_{\d}g\right\|_p\lesssim_p\delta^{\frac{1}{p}}\left\|g\right\|_p$
and $\left\|R*g\right\|_p\lesssim_p\left\|g\right\|_p$, we conclude
that \eqref{tree} holds.
\end{proof}
\end{l1}

\begin{d0}\label{sep}
Fix a number $\d\in(0,1)$. Let $\p_1$ and $\p_2$ be two trees with
tops $P_1=[\a_1^1,\a_1^2, \ldots \a_1^d,I_1]$ and respectively
$P_2=[\a_2^1,\a_2^2,\ldots \a_2^d,I_2]$; we say that $\p_1$ and
$\p_2$ are \emph{($\d$-)separated} if $\newline$ either $I_1\cap
I_2=\es$ or else$\newline$i) $P=[\a^1,\a^2,\ldots
\a^d,I]\in\p_1\:\:\&\:\:I\subseteq
I_2\:\:\:\:\Rightarrow\:\:\:\left\lceil
\Delta(P,P_2)\right\rceil<\d\:,$ $\newline$ii) $P=[\a^1,\a^2,\ldots
\a^d,I]\in\p_2\:\:\&\:\:I\subseteq
I_1\:\:\:\:\Rightarrow\:\:\:\left\lceil
\Delta(P,P_1)\right\rceil<\d\:.$
\end{d0}

\noindent{\bf Notation:} Let $\p_1$ and $\p_2$ be two trees as in
Definition \ref{sep}. Take $q_j$ to be the central polynomial of
$P_j$ ($j\in\{1,2\}$), set $q_{1,2}=q_1-q_2$, and then define

 \begin{itemize}
\item $I_s$ - the {\bf separation set} (relative to the intersection) of $\p_1$ and $\p_2$ by
$$I_s=\I_s\left(\eta_{1,2},c(d)\d^{-1},q_{1,2},\tilde{I}_1\cap\tilde{I}_2\right)\:,$$
where $\eta_{1,2}=\eta(\tilde{I}_1\cap\tilde{I}_2)$;
\item $I_c$ - the {\bf ($\ep$-)critical
intersection set} (between $\p_1$ and $\p_2$) by
$$I_c=\bigcup_{j=1}^{r}\I_c\left(\eta(I_s^j),w(I_s^j),q_{1,2},I_s^j\right)\:,$$ where $\ep$ is some small fixed positive real
number and $I_s=\bigcup_{j=1}^{r}I_s^j$ is the decomposition of
$I_s$ into maximal disjoint intervals $\{I_s^j\}_{j\in\{1,\ldots
r\}}$ with $r\leq 2d$.
  \end{itemize}
\noindent{\bf Observation 5.}  It is important to notice the
following three properties of our above-defined sets; these
facilitate the adaptation of the reasonings involved in the proofs
of Lemmas 2 and 4 to those of the corresponding lemmas in \cite{q}:

1) for all dyadic $J\subset \tilde{I}_1\cap\tilde{I}_2$ such that
$I_s\cap 5J=\emptyset$ we have $(c(d)\leq d^d)$
$$\inf_{x\in J}|q_{1,2}(x)|\leq \sup_{x\in J}|q_{1,2}(x)|\leq c(d)\inf_{x\in J}|q_{1,2}(x)|\:.$$
\indent 2) $\forall\:\:P\in\p_1\cup \p_2$ and $j\in\{1,\ldots r\}$ if
$I_s^j\cap 5\tilde{I}_P\not=\emptyset$ then $|I_P|>|I_s^j|$.

3) $\forall\:\:P\in\p_1\cup \p_2$ we have (for $\ep$ properly
chosen) $|\tilde{I}_P\cap I_c|<{\d}^{\frac{1}{100d}}|I_P|$.

Indeed, these facts are an easy consequence of the results mentioned
in the Appendix and the way in which $I_s$ and $I_c$ are defined.

(Remark that property 1) above implies the following relation:
\begin{quote}
$\:\:\forall\:P\in\p_1$ such that $I_s\cap
5\tilde{I}_P=\emptyset$ and $I_P\subset I_2$ we have
$$\textrm{Graph}(q_2)\cap\left(c(d)\d^{-1}\right)\widehat{P}=\emptyset\:.$$
\end{quote}
Of course, the same is true for the symmetric relation, {\it i.e.}
replacing the index $1$ with $2$ and vice versa.)

\begin{l1}\label{sept}
Let $\left\{\p_j\right\}_{j\in\left\{1,2\right\}}$  be two separated
trees with tops $\\ P_j=[\a^1_j,\a^2_j,\ldots \a^d_j,I_0]$. Then,
for any $f,\:g\in L^{2}(\TT)$ and $n\in \N$, we have that
\beq\label{v21} \left|\left\langle
{T^{\p_1}}^*f,\,{T^{\p_2}}^*g\right\rangle\right|\lesssim_{n,d}{\d}^n\left\|f\right\|_{L^{2}(\tilde{I}_0)}\left\|g\right\|_{L^{2}(\tilde{I}_0)}+\left\|\chi_{I_c}{T^{\p_1}}^*f\right\|_2\left\|\chi_{I_c}{T^{\p_2}}^*g\right\|_2\:.
\eeq
\end{l1}
\begin{proof}
In what follows we intend to adapt the methods described in the
proof of Lemma 2 of \cite{q} to our context. For this, we need first
to modify the definition of the sets $\left\{A_l\right\}_l$; more
exactly we follow the procedure below: Let
$I_s=\bigcup_{j=1}^{r}I_s^j$ be the decomposition of $I_s$ into
maximal disjoint intervals ($r\leq 2d$); without loss of generality
we may suppose\footnote{We made here the convention
$I_s^0=I_s^{r+1}=\emptyset$} that $\left\{I_s^j\right\}_j$ are
placed in consecutive order with $I_s^{j+1}$ located to the right of
$I_s^j$. For a fixed $j\in \left\{0,\ldots r\right\}$, let $W_j$ be
the standard Whitney decomposition of the set
$[0,1]\cap\left(I_s^j\right)^c\cap\left(I_s^{j+1}\right)^c$ with
respect to the set $I_s^j\cup I_s^{j+1}$; we take $\tilde{W}_j$ to
be the ``large scale" version of $W_j$, which is obtained as
follows:

- take the union of all the intervals in $W_j$ of length strictly smaller than $c(d)|I_s^j|$ that approach $I_s^j$ and
denote it by $R_j$ (we can do this in such a way that $R_j$ can be written as a union of at most two dyadic intervals, each one of length $c(d)|I_s^j|$ );

- apply the same procedure to obtain $R_{j+1}$;

- the rest of the intervals belonging to $W_j$ remain unchanged and are transferred to $\tilde{W}_j$.

Define $\W_j:=\tilde{W}_j\cup I_s^j\cup I_s^{j+1}$ and observe that this is a partition of $[0,1]$.

Finally, we take $\W$ to be the common refinement of the partitions
$\W_j$, $j\in\{0,\ldots r\}$. Take now $A_0=\bigcup_{I\in\W\atop{
\bar{I}\cap\bar{I_s}\neq\emptyset}} I$ and set
$$\W=A_0\cup \bigcup_{l=1}^k A_l\:.$$

Then, for $l\in\{1,\ldots k\}$ and $m\in\{1,2\} $ we define the sets
$$S_{m,l}:=\left\{P\in\p_m\:|\:I_P\subset A_l\:\:\&\:\:|I_P|\leq\frac{|A_l|}{20}\right\}\:.$$

Also, we take
$S_{m,0}:=\p_m\setminus\left(\bigcup_{l=1}^{k}S_{m,l}\right)\:.$

With this done, for $l\in\{0,\ldots k\}$, we set
$$T^{*}_{m,l}=\sum_{P\in S_{m,l}}T^{*}_{P}$$
and deduce that
$$\left\langle{T^{\p_1}}^*,\:{T^{\p_2}}^*
\right\rangle=\sum_{n,l=0}^{k}\left\langle{T_{1,l}}^*,\:{T_{2,n}}^*
\right\rangle\:. $$

Now, as intended, we may follow the same steps as in \cite{q}, Lemma
2 and show that \beq\label{nonc}
\sum_{l=0}^k\sum_{n=1}^k\left|\left\langle
{T_{1,l}}^*f,\,{T_{2,n}}^*g\right\rangle\right|\lesssim_{n,d}
{\d}^n\left\|f\right\|_{L^{2}(\tilde{I}_0)}\left\|g\right\|_{L^{2}(\tilde{I}_0)}\:;\eeq
\beq\label{crit}
\left|\left\langle{T_{1,0}}^*f,\,{T_{2,0}}^*g\right\rangle\right|\lesssim_{n,d}
{\d}^n\left\|f\right\|_{L^{2}(\tilde{I}_0)}\left\|g\right\|_{L^{2}(\tilde{I}_0)}+
\left\|\chi_{I_c}{T^{\p_1}}^*f\right\|_2\left\|\chi_{I_c}{T^{\p_2}}^*g\right\|_2\:.\eeq
 finishing our proof.
\end{proof}

\begin{d0}\label{nor}
A tree $\p$ with top $P_0=[\a_0^1,\a_0^2,\ldots
\a_0^d,I_0]$ is called \emph{normal} if $\newline$
$P=[\a^1,\a^2,\ldots
\a^d,I]\in\p\:\:\:\:\:\:\Rightarrow\:\:\:\:\:\:|I|\leq\frac{\d^{100}}{K^{2M}}|I_0|\:\:\:\&\:\:\:\:dist(I,\partial
I_0)>20\frac{\d^{100}}{K^{2M}}|I_0|$.
\newline(Here $K,\:M\in\N$ are some fixed large constants and $ \partial I_0$ designates the
boundary of $I_0$.)
\end{d0}
\noindent{\bf Observation 6.} Notice that if $\p$ is a normal tree
then $$supp\:{{T^{\p}}^{*}f}\subseteq \left\{x\in
I_0\:|\:dist(x,\partial
I_0)>10\frac{\d^{100}}{K^{2M}}|I_0|\right\}\:.$$ Also, it is worth
mentioning that the constant $M$ introduced above is exactly that
appearing in the variational estimates used in Section 6.
\begin{d0}\label{row}
A \emph{row} is a collection $\p=\bigcup_{j\in \N}\p^{j}$ of normal
trees $\p^{j}$ with tops
$P_{j}=[\a^1_{j},\a^2_j,\ldots \a^d_j,I_j]$ such that the
$\left\{I_j\right\}$ are pairwise disjoint.
\end{d0}

The proofs of the next two lemmas require no nontrivial
modifications from the corresponding proofs in \cite{q}.
\begin{l1}\label{rt}
Let $\p$ be a row as above, let $\p'$ be a tree with top $\newline$
$P'=[{\a^1_{0}}',{\a^2_{0}}',\ldots {\a^d_{0}}',I_{0}']$ and suppose
that $\forall\: j\in \N$, $I_0^{j}\subseteq I_0'$ and $\p^{j},\:\p'$
are separated trees; denote by $I_c^{j}$ the critical intersection
set between each $\p^{j}$ and $\p'$.

Then for any $f,\:g\in L^{2}(\TT)$ and $n\in \N$ we have that
$$\left|\left\langle {T^{\p'}}^*f,{T^{\p}}^*g\right\rangle\right|
\lesssim_{n,d}{\d}^n\left\|f\right\|_{2}\left\|g\right\|_{2}+
\left\|\sum_{j}\chi_{I_c^{j}}{T^{\p'}}^*f\right\|_2\left\|\sum_{j}\chi_{I_c^{j}}{T^{\p^{j}}}^*g\right\|_2\:.$$
\end{l1}

\begin{l1}\label{ed}
Let $\p\:$ be a tree with top
$P_{0}=[\a^1_{0},\a^2_{0},\ldots \a^d_{0},I_{0}]$; suppose also that
we have a set $A\subseteq \tilde{I_0}$ with the property that \beq
\label{masx} \:\exists\:\d\in(0,1)\:st\:\:\:\:\:\:\:\forall
\:P=[\a^1,\a^2,\ldots
\a^d,I]\in\p\:\:\operatorname{we\:have\:}\:\:|I^{*}\cap A|\leq \d
|I|. \eeq

Then  $\forall\:f\in L^2(\TT)$ we have \beq \label{cut}
\left\|\chi_{A}{T^{\p}}^{*}f\right\|_2\lesssim_{d}\d^{\frac{1}{2}}\left\|f\right\|_{2}.
\eeq
\end{l1}

$\newline${\bf Proof of Proposition 2}

$\newline$

As in \cite{q}, we can ignore the behavior of our operator on the
set $F=\bigcup_{j}\left\{x\in I_j\:|\operatorname{dist}(x,\partial
I_j))\leq100\frac{\d^{100}}{K^{2M}}|I_j|
\right\}=^{def}\bigcup_{j}F_j\:,$ since
$$|F|\leq\sum_{j}|F_j|\lesssim\sum_j |I_j|\frac{\d^{100}}{K^{2M}}\lesssim\frac{\d^{50}}{K^{M}}\:.$$

Now on $F^{c}$ we will use the estimates obtained in Lemma 3. But
before this, we first need to remove a few tiles\footnote{If for a
tree $\p_j$ there are multiple tiles with the same time interval,
then we will take their (geometric) union and consider it as a
single tile.} from each tree $\p_j$.

For $\p=\bigcup_{j}\p_j$ and $L=\log{(K^{100M}\d^{-100M})}\:$ we let
$$\p^{+}:=\left\{P\in\p\:|\:\operatorname{there\:is\:no\:chain\:}
 P<P_1<...<P_{L}\:\operatorname{with\:all\:}P_j\in\p\right\}$$
and
$$\p^{-}:=\left\{P\in\p\:|\:\operatorname{there\:is\:no\:chain\:}
 P_1<P_2<...<P_{L}<P\:\operatorname{with\:all\:}P_j\in\p\right\}\:.$$

Now, it is easy to see that each such set can be split into at most
$L$ subsets with no two comparable tiles inside the same subset.
Consequently, using Proposition 1, we deduce that\footnote{$\eta$
may change from line to line.}
$$\left\|T^{\p^{+}}\right\|_p\:,\:\left\|T^{\p^{-}}\right\|_p\lesssim_{p,d} L{\d}^{\eta(1-\frac{1}{p^*})}\lesssim{\d}^{\eta(1-\frac{1}{p^*})}\log{K}\:.$$

We remove all the above-mentioned sets from our collection $\p$ and
decompose this new set as follows:
$$\p=\bigcup_{j}\p^{0}_j\:\:\operatorname{where}\:\: \p^{0}_j=\p_j\cap\p\:.$$

This new collection $\p$ has the following properties:
$\newline\:1)\:\forall\:P\in\p^{0}_j$,
$\:|I_P|\leq\frac{\d^{100M}}{K^{100M}}|I_j|$
$\newline\:2)\:\forall\:j\not=k$, the trees $\p^{0}_j$ and
$\p^{0}_k$ are $\d'$-separated where
$\d'=\frac{\d^{100M}}{K^{100M}}\:.$

Splitting now each $\p^{0}_j=\p^{N}_j\cup\p^{C}_j$, with
$$\p^{C}_j=^{def}\left\{P\in\p^{0}_j\:|\:I_P\subseteq F_j\right\}\:,$$
we conclude that $\:\left\{\p^{N}_j\right\}_{j}$ represents a
collection of normal, $\d'$-separated trees, while for the remaining
parts of the trees we have the relation
$$\operatorname{supp}\:T^{\p^{C}_j}\subset F_j\:.$$

 Consequently, on $F^c$ we have that
$$T^{\p}f\:=\:\sum_{j}T^{\p^{N}_j}f\:,$$
and so our conclusion reduces to \beq\label{enough}
\left\|\sum_{j}T^{\p^{N}_j}f\right\|_p\lesssim_{p,d}
\left(K^{\frac{1}{p^*}-\frac{1}{p'}}\:\d^{\frac{1}{p}-(1+\rho)(\frac{1}{p^*}-\frac{1}{p'})}\right)\left\|f\right\|_p\:.
 \eeq
 Now, as in \cite{q}, we may divide
$\bigcup_{j}\p^{N}_j$ into a union of at most $S:=K\d^{-(1+\rho)}$
rows, $\r_1\:,\:\r_2\:,...\r_S$. Using the same techniques as in
\cite{2} and further \cite{q}, one can show (applying the
Cotlar-Stein Lemma together with Lemmas 1 - 4) that \beq\label{L2}
\left\|\sum_{j}{T^{\p^{N}_j}}^{*}f\right\|_2=\left\|\sum_{j=1}^S
{T^{{\r}_j}}^{*}f\right\|_2\lesssim_{d}
\left\{\sum_{j=1}^S\left\|{T^{{\r}_j}}^{*}f\right\|_2^2\right\}^{1/2}\:\left(\lesssim_d\d^{\frac{1}{2}}\left\|f\right\|_2\right)\:.
\eeq

 On the other hand, from the triangle inequality, we trivially have
\beq\label{L1} \left\|\sum_{j=1}^S {T^{{\r}_j}}^{*}f\right\|_1\leq
\sum_{j=1}^S\left\|{T^{{\r}_j}}^{*}f\right\|_1
\:\:\:\:\:\:\:\:\:\:\: \left\|\sum_{j=1}^S
{T^{{\r}_j}}^{*}f\right\|_{\infty}\leq
\sum_{j=1}^S\left\|{T^{{\r}_j}}^{*}f\right\|_{\infty}\:. \eeq

Now, using interpolation between \eqref{L2} and \eqref{L1}, we
obtain \beq\label{Lp} \left\|\sum_{j=1}^S
{T^{{\r}_j}}^{*}f\right\|_{p'}\lesssim_{p,d}
\left\{\sum_{j=1}^S\left\|{T^{{\r}_j}}^{*}f\right\|_{p'}^{p^*}\right\}^{\frac{1}{p^*}}\:.
\eeq

If we denote $E_j:=\textrm{supp}\:T^{{\r}_j}$ then, by Lemma 1, we
deduce
$$\left\|{T^{{\r}_j}}^{*}f\right\|_{p'}\lesssim_{p',d}
\d^{\frac{1}{p}}\left(\int_{E_j}|f|^{p'}\right)^{\frac{1}{p'}}\:.$$

Now, taking into account the fact that $\{E_j\}_j$ are disjoint, we
conclude (using H\"older)
$$\left\{\sum_{j=1}^S\left(\int_{E_j}|f|^{p'}\right)^{\frac{p^*}{p'}}\right\}^{\frac{1}{p^*}}\lesssim_{p,d}
S^{\frac{1}{p^*}-\frac{1}{p'}}\left(\int|f|^{p'}\right)^{\frac{1}{p'}}\:,$$
which ends our proof.

 \begin{flushright}
$\Box$
\end{flushright}

\section{\bf Appendix - Results on the $L^{\infty}-$distribution of polynomials}

{\bf Lemma A.} {\it If $q\in\Q_{d-1}$ and $I,\:J$ are some (not
necessarily dyadic) intervals obeying $I\supseteq J$, then there
exists a constant $c(d)$ such that
$$ \left\|q\right\|_{L^{\infty}(I)}\leq
c(d)\left(\frac{|I|}{|J|}\right)^{d-1}\left\|q\right\|_{L^{\infty}(J)}\:.$$}
$\newline${\it Proof.} Let $\{x_J^k\}_{k\in\{1,\ldots d\}}$ be
obtained as in the
 procedure described in Section 2. Then, since $q\in\Q_{d-1}$, for any $x\in I$ we have that
 $$q(x):=\sum_{j=1}^{d}\frac{\prod_{k=1\atop{k\not=j}}^{d}(x-x_J^k)}{\prod_{k=1\atop{k\not=j}}^{d}(x_J^j-x_J^k)}\:q(x_J^j)\:.$$
As a consequence,
$$\left\|q\right\|_{L^{\infty}(I)}\leq
d\:\left\|q\right\|_{L^{\infty}(J)}\sup_{j\atop{x\in I}}\left|
\frac{\prod_{k=1\atop{k\not=j}}^{d}(x-x_J^k)}{\prod_{k=1\atop{k\not=j}}^{d}(x_J^j-x_J^k)}\right|\leq
d\:\left\|q\right\|_{L^{\infty}(J)}\frac{|I|^{d-1}}{(|J|/d)^{d-1}}\:.$$

\begin{flushright}
$\Box$
\end{flushright}

{\bf Lemma B.} {\it If $q\in\Q_{d-1},\:\eta>0$ and $I\subset\TT$
some (dyadic) interval, then \beq\label{levels} |\{y\in
I\:|\:|q(y)|<\eta\}|\leq c(d)
\left(\frac{\eta}{\left\|q\right\|_{L^{\infty}(I)}}\right)^{\frac{1}{d-1}}|I|\:.
\eeq } $\newline${\it Proof.} The set $A_{\eta}=\{y\in
I\:|\:|q(y)|<\eta\}$ is the pre-image of $(-\eta,\eta)$ (an open
set) under a polynomial of degree $d-1$, so it can be written as
$$A_{\eta}=\bigcup_{k=1}^r J_k(\eta)\:,$$
where $r\in\N,\:r\leq d-1$ and $\{J_k(\eta)\}_k$ are open intervals.
Now all that remains is to apply the previous lemma with
$J=J_k(\eta)$ for each $k$. \hfill $\Box$

{\bf Lemma C.} {\it If $P=[\a^1,\a^2,\ldots \a^d,I]\in\P$ and $q\in
P$, then
$$\left\|q-q_P\right\|_{L^{\infty}(\tilde{I})}\leq
c(d)\:|I|^{-1}\:.$$} $\newline${\it Proof.}
Set $u:=q-q_P$; then, since both $q,\:q_P\in P$, we deduce (for all
$k\in\{1,\ldots d\}$):
$$u(x_I^k)\in [-|I|^{-1},|I|^{-1}]\:.$$
On the other hand,
$$u(x):=\sum_{j=1}^{d}\frac{\prod_{k=1\atop{k\not=j}}^{d}(x-x_I^k)}{\prod_{k=1\atop{k\not=j}}^{d}(x_I^j-x_I^k)}\:u(x_I^j)\:\:\:\:\:\:\:\:\:\forall\:\:x\in I\:.$$
Then, proceeding as in Lemma A, we conclude
$$\left\|u\right\|_{L^{\infty}(I)}\leq d\:
|I|^{-1}\frac{|I|^{d-1}}{(|I|/d)^{d-1}}\leq {d}^d\: |I|^{-1}\:.$$

\begin{flushright}
$\Box$
\end{flushright}

\end{document}